\documentclass[10pt]{amsart}
\usepackage{amssymb}

\usepackage{amscd}
\usepackage{amsthm}
\usepackage[dvips]{graphicx}

\usepackage[all]{xy}

\newtheorem{Thm}{Theorem}[section]
\newtheorem{Lem}[Thm]{Lemma}
\newtheorem{Def}[Thm]{Definition}
\newtheorem{Cor}[Thm]{Corollary}
\newtheorem{Prop}[Thm]{Proposition}
\newtheorem{Ex}[Thm]{Example}
\newtheorem{Rem}[Thm]{Remark}

\newtheorem{Crit}[Thm]{Criterion}

\title[Batalin-Vilkovisky]{
The Hochschild cohomology ring  of a Frobenius algebra with
semisimple Nakayama automorphism is a   Batalin-Vilkovisky
algebra}

\author{Thierry Lambre}

\author{Guodong Zhou}

\author{Alexander Zimmermann}

\address{Thierry Lambre\\
 \newline Laboratoire de Math\'ematiques
 \newline UMR 6620 du  CNRS
\newline Universit\'e Blaise Pascal
\newline  63177 Aubi\`ere Cedex
\newline France}
  \email{thierry.lambre@math.univ-bpclermont.fr}

\address{Guodong Zhou
 \newline Department of Mathematics
 \newline Shanghai Key laboratory of PMMP
\newline East China Normal University
\newline  Dong Chuan Road 500
\newline Shanghai 200241
 \newline P.R.China}
  \email{gdzhou@math.ecnu.edu.cn}

  \address{Alexander Zimmermann
\newline D\'epartement de Math\'ematiques et LAMFA (UMR 7352 du CNRS),
\newline Universit\'e de Picardie,
\newline 33 rue St Leu,
\newline F-80039 Amiens Cedex 1,
\newline France}
\email{alexander.zimmermann@u-picardie.fr}

\date{version of \today}

\newcommand{\comp}{\mathop{\raisebox{+.3ex}{\hbox{$\scriptstyle\circ$}}}}
\newcommand{\xto}{\xrightarrow}
\newenvironment{Proof}[1][Proof]{\begin{trivlist}
\item[\hskip \labelsep {\bfseries #1}]}{\flushright
$\Box$\end{trivlist}}

\newcommand{\lra}{\longrightarrow}

\newcommand{\ra}{\rightarrow}
\newcommand{\sdp}{\times\kern-.2em\vrule height1.1ex depth-.05ex}
\newcommand{\epi}{\lra \kern-.8em\ra}

\newcommand{\N}{{\mathbb N}}

\newcommand{\ol}{\overline}

\newcommand{\cal}{\mathcal}

\newcommand{\ot}{\otimes}

 \normalbaselines
\begin{document}
\renewcommand{\thefootnote}{\alph{footnote}}
\renewcommand{\thefootnote}{\alph{footnote}}
\setcounter{footnote}{-1} 
\renewcommand{\thefootnote}{\alph{footnote}}
\setcounter{footnote}{-1} 
\begin{abstract} Analogous to a recent result of N. Kowalzig and
 U. Kr\"{a}hmer for twisted Calabi-Yau algebras,  we show that the Hochschild cohomology ring
 of a Frobenius algebra with semisimple Nakayama automorphism is a Batalin-Vilkovisky algebra,
 thus generalizing a result of T.Tradler for finite dimensional symmetric algebras.  We
 give a criterion to determine when a Frobenius algebra given by quiver with relations has
 semisimple Nakayama automorphism and apply it to some known classes of tame Frobenius
 algebras. We also provide ample examples   including quantum complete intersections,
 finite dimensional Hopf algebras defined over an algebraically closed field of
 characteristic zero and Koszul duals of Koszul Artin-Schelter regular algebras of dimension three.
\end{abstract}

\maketitle

%\tableofcontents

\section*{Introduction}

Let $A$ be an associative algebra over a field $k$. The Hochschild
cohomology groups $HH^{*}(A)$ of $A$ has a very rich structure. It  is a graded
commutative algebra
via the cup product or the Yoneda product, and  it has a graded  Lie bracket
of degree $-1$ so that it becomes  a graded Lie
algebra; these make ${HH}^*(A)$ a Gerstenhaber
algebra \cite{Gerstenhaber}. Furthermore,  Hochschild homology groups are
endowed with two actions by the Hochschild cohomology algebra, which make the
Hochschild homology groups to be graded modules and graded Lie modules over the
Hochschild cohomology algebra.  These structures are summarized by the notion of
a differential calculus (see \cite{GelfandDaleskiiTsygan} and \cite{TamarkinTsygan});
we explain in  detail these structures in the first section.

During several decades,  a new structure in Hochschild theory has been
extensively studied in topology and mathematical physics, and recently this was
introduced into algebra, the so-called Batalin-Vilkovisky structure. Roughly
speaking a Batalin-Vilkovisky (aka. BV) structure is an operator on Hochschild
cohomology which squares to zero and which, together with the cup product, can
express the Lie bracket.  A BV structure only exists over Hochschild cohomology
of certain special classes of algebras.

T.~Tradler first found  that the Hochschild cohomology algebra of a
finite dimensional symmetric algebra, such as a group algebra of a finite group,
is a BV algebra \cite{Tradler}; for later proofs,  see e.g.
 \cite{EuSchedler}\cite{Menichi}.
For a Calabi-Yau algebra, V.~Ginzburg showed in \cite{Ginzburg} that the Hochschild
cohomology of a Calabi-Yau algebra also has a Batalin-Vilkovisky
structure.

Inspired by the result of V. Ginzburg, the first author introduced in \cite{Lambre}
the notion of a differential calculus with duality. This notion explains when BV
structure exists and  unifies the two known cases of symmetric algebras and
Calabi-Yau algebras. Recently as an application of this notion, N.~Kowalzig
and U.~Kr\"{a}hmer  \cite[Theorem 1.7]{KowalzigKrahmer} proved that the
Hochschild cohomology ring of a twisted Calabi-Yau algebra is also a
Batalin-Vilkovisky algebra, provided a certain algebra automorphism is semisimple.

The main result of this paper is an
analogous statement for Frobenius algebra with semisimple Nakayama
automorphism. Our main result reads as follows.

\begin{Thm}\label{Thm:MainResult} Let $A$ be a Frobenius algebra with semisimple
Nakayama automorphism. Then  the Hochschild cohomology ring $HH^*(A)$ of
$A$ is  a Batalin-Vilkovisky algebra.
\end{Thm}

Observe that the semisimplicity is an open condition, and that  any finite
dimensional self-injective algebra defined over an algebraically closed field
is Morita equivalent to its basic algebra which is   a Frobenius algebra.
Hence our main result shows that  the Hochschild cohomology rings of most
self-injective algebras are BV algebras.

The paper is organised as follows.  In  Section~\ref{Sectionone}  we explain
the formalism  of
Tamarkin-Tsygan calculi, calculi with duality and Batalin-Vilkovisky structures.
Section~\ref{Sectiontwo}
develops the Tamarkin-Tsygan structure on the Hochschild cohomology associated
with an automorphism of an algebra. We show that when the   Nakayama
automorphism of a Frobenius algebra is  diagonalisable, then there is a
differential calculus with duality which is a key ingredient of the proof
of our main result.
Section~\ref{Sectionthree} then studies the special case of a Frobenius algebra.
We provide a proof of the main result in Section~\ref{Sectionfour}.
Section~\ref{Sectionfive} contains many examples of Frobenius algebras
with semisimple Nakayama automorphisms. For  a Frobenius algebra given
in terms of quiver with relations, we give a very useful combinatorial
criterion to guarantee the  semisimplicity of the Nakayama automorphism
and apply it to some classes of tame Frobenius algebras.  We then include
other examples      such as quantum complete intersections, finite dimensional
Hopf algebras and Koszul duals of Artin-Schelter regular algebras.

Throughout this paper, $\ot$ is an abbreviation for $\ot_k$ for $k$ being a chosen  base field.

\begin{Rem}\rm
After having finished this paper we learned that independently Y.~Volkov proved in~\cite{Volkov}
a similar result with completely different methods.
He works directly over Hochschild cohomology by defining an operator
analogous to Tradler's
operator twisted by the Nakayama automorphism. However,  our method uses
the concept of Tamarkin-Tsygan calculi.

We are grateful to Y.~Volkov for pointing out an error in a previous version of this paper.

The last two authors are partially supported by the exchange program
STIC-Asie 'ESCAP' financed by the French Ministry of Foreign Affairs.
The second author is supported    by Shanghai Pujiang
Program (No.13PJ1402800),  by National Natural Science Foundation of
China (No.11301186) and by the Doctoral Fund of Youth Scholars of
Ministry of Education of China (No.20130076120001).
\end{Rem}

\section{Tamarkin-Tsygan calculus, duality and Batalin-Vilkovisky structure}

\label{Sectionone}

\subsection{Gerstenhaber Algebras}

First  we recall the definition of Gerstenhaber algebras and of
differential calculi.

\begin{Def}\label{Def:Gerstenhaber}
A {\em Gerstenhaber algebra} over a field $k$ is the data $(\mathcal{H}^*, \cup, [\ ,\ ])$,  where  $\mathcal{H}^*=\oplus_{n\in \mathbb{Z}
}\mathcal{H}^n$ is  a   graded
$k$-vector space equipped with two bilinear maps: a cup product of
degree zero
$$\cup\colon \mathcal{H}^n \times \mathcal{H}^m \rightarrow \mathcal{H}^{n+m}, \quad (\alpha,\, \beta) \mapsto
\alpha\cup \beta$$ and a Lie bracket of degree $-1$
$$[\ ,\ ]\colon \mathcal{H}^n\times \mathcal{H}^m\rightarrow \mathcal{H}^{n+m-1}, \quad (\alpha,\, \beta)\mapsto
[\alpha,\, \beta]$$ such that
\begin{itemize}
\item[(i)] $(\mathcal{H}^*,\,  \cup)$ is a graded commutative
associative algebra with unit $1\in \mathcal{H}^0$, in particular, $\alpha\cup
\beta=(-1)^{|\alpha||\beta|}\beta\cup \alpha$;

\item[(ii)] $(\mathcal{H}^*[-1],\, [\ ,\,\ ])$ is a graded Lie algebra, that is,
$$[\alpha,\, \beta]=-(-1)^{(|\alpha|-1)(|\beta|-1)}[\beta,\, \alpha]$$ and
$$(-1)^{(|\alpha|-1)(|\gamma|-1)}[[\alpha,\, \beta],\, \gamma]+(-1)^{(|\beta|-1)(|\alpha|-1)}[[\beta,\, \gamma],\, \alpha]
+(-1)^{(|\gamma|-1)(|\beta|-1)}[[\gamma,\, \alpha],\, \beta]=0; $$

\item[(iii)] for each $\alpha\in{\mathcal H}^*[-1]$ the map $[\alpha,-]$ is a graded derivation of the algebra $(\mathcal H^*,\cup)$, or more precisely
$$[\alpha, \, \beta \cup \gamma]=[\alpha,\, \beta]\cup \gamma+(-1)^{(|\alpha|-1)|\beta|}\beta\cup [\alpha,\, \gamma],$$

\end{itemize}
where  $\alpha, \beta, \gamma$ are arbitrary homogeneous elements
in $ \mathcal{H}^*$ and $|\alpha|$ is the degree of the homogeneous element $\alpha$.
\end{Def}

\begin{Rem}\rm \label{Rem:ExtensionDefinitionOfGertstenhaberAlgebras} Let $k'$
be a field extension of $k$. Then for a Gerstenhaber algebra $\mathcal{H}^*$
over $k$, $\mathcal{H}^*\ot k'$ is a Gerstenhaber algebra over $k'$. In fact,
for homogenous elements $\alpha, \beta\in \mathcal{H}^*$ and $\lambda, \mu\in k'$, define
$$(\alpha\ot \lambda) \cup (\beta\ot \mu)=(\alpha\cup\beta)\ot (\lambda\mu)$$
and $$[\alpha\ot \lambda,\, \beta\ot \mu]=[\alpha, \beta]\ot (\lambda\mu).$$
These two operations endow $\mathcal{H}^*\ot k'$ with  a Gerstenhaber algebra structure over $k'$.
\end{Rem}

\subsection{Tamarkin-Tsygan Calculi}

\begin{Def} \label{Def:Calculus} A differential calculus or a Tamarkin-Tsygan calculus is the data
$(\mathcal{H}^*, \cup,   [\ ,\ ],\mathcal{H}_*,  \cap, B)$
of $\mathbb{Z}$-graded vector spaces satisfying the following
properties:

\begin{itemize}

\item[(i)] $(\mathcal{H}^*, \cup, [\ ,\ ])$ is a Gerstenhaber
algebra;

\item[(ii)] $\mathcal{H}_*$ is a graded module over $(\mathcal{H}^*,
\cup)$ via the map $\cap:   \mathcal{H}_r\ot \mathcal{H}^p\to
\mathcal{H}_{r-p},    z\ot \alpha \mapsto   z\cap \alpha$ for $z\in
\mathcal{H}_r$ and $\alpha\in \mathcal{H}^p$. That is, if we
denote $\iota_\alpha(z)= (-1)^{rp} z\cap\alpha$, then
$\iota_{\alpha\cup \beta}=\iota_{\alpha}\iota_{\beta}$;

\item[(iii)] There is a map $B: \mathcal{H}_*\to
\mathcal{H}_{*+1}$ such that $B^2=0$ and we have the Tamarkin-Tsygan relation
$$[L_\alpha, \ \iota_{\beta}]_{gr}=\iota_{[\alpha, \beta]}$$
where we denote
$$L_\alpha=[B, \iota_\alpha]_{gr}=B
\iota_{\alpha}-(-1)^{|\alpha|} \iota_{\alpha} B.$$

\end{itemize}

\end{Def}

One of the first examples of differential calculi is  Hochschild
theory.

The cohomology theory of associative algebras was introduced by G.
Hochschild (\cite{Hochschild}).   Given a $k$-algebra $A$, its
Hochschild cohomology groups of $A$ with coefficients in a
bimodule $M$ are defined as $HH^n(A, M)= \mathrm{Ext}^n_{A^e}(A,
M)$ for $n \geq 0$, where $A^e=A \otimes A^{\mathrm{op}}$ is the
enveloping algebra of $A$, and the Hochschild homology groups  of
$A$ with coefficients in   $M$ are defined to be $HH_n(A,
M)=Tor^{A^e}_n (A, M)$ for $n\geq 0$. We shall write
$HH_n(A)=HH_n(A, A)$ and $HH^n(A)=HH^n(A, A)$.

Since $A$ is unitary, denote by $1_A$ its unity and write $\ol{A}=A/(k\cdot 1_A)$.
For $a\in A$, write $\ol{a}$ for its image in $\ol{A}$. There is a projective
resolution of $A$ as an $A^e$-module, the so-called   \textit{normalized bar resolution}
$\mathrm{Bar}_*(A)$, whose $r$-th term
is given by $\mathrm{Bar}_r(A)=A\ot \ol{A}^{\ot r}\ot A$ for $p\geq 0$  and for which
the differential $b'_r: \mathrm{Bar}_r(A)\to \mathrm{Bar}_{r-1}(A)$ sends
$a_0\otimes \ol{a_1}\otimes\cdots\otimes\ol{a_r}\ot  a_{r+1}$ to
$$\begin{array}{l}
a_0a_1\ot \ol{a_2}\ot \cdots \ot \ol{a_r}\ot a_{r+1}+\sum_{i=1}^{r-1}(-1)^i
a_0 \otimes \ol{a_1}\otimes \cdots\otimes \ol{a_{i-1}}\otimes  \ol{a_ia_{i+1}}\otimes
\ol{a_{i+2}}\otimes \cdots\otimes \ol{a_r}\otimes  a_{r+1}\\+(-1)^r a_0\ot \ol{a_1} \ot \cdots \ot \ol{a_{r-1}}\ot a_ra_{r+1} \end{array}
$$
 for all $a_0, \cdots, a_{r+1} \in A$.

The complex which is used to compute the Hochschild cohomology is
$C^*(A, M)=\mathrm{Hom}_{A^e}(\mathrm{Bar}_*(A), M)$. Note that
for each $r\geq 0$, $$C^r(A, M)=\mathrm{Hom}_{A^e}(A\ot \ol{A}^{\otimes
r}\ot A, M)\simeq\mathrm{Hom}_k(\ol{A}^{\otimes r}, M).$$ We identify
$C^0(A, M)$ with $M$.  Thus $C^*(A, M)$ has the following form:
$$C^*(A, M) \colon M \xto{b^0} \mathrm{Hom}_k(\ol{A}, M) \to \cdots \to \mathrm{Hom}_k(\ol{A}^{\otimes r}, M)
\xto{b^r} \mathrm{Hom}_k(\ol{A}^{\otimes (r+1)}, M)\to \cdots .$$
Given $f$ in $\mathrm{Hom}_k(\ol{A}^{\otimes r}, M)$, the map
$b^r(f)$ is defined by sending $\ol{a_1}\otimes \cdots\otimes
\ol{a_{r+1}}$ to
\[
(-1)^{r+1}a_1 \cdot f(\ol{a_2}\otimes \cdots\otimes
\ol{a_{r+1}})+\sum_{i=1}^r(-1)^{r+1-i} f(\ol{a_1}\otimes \cdots \otimes
\ol{a_{i-1}}\otimes \ol{a_ia_{i+1}} \otimes \ol{a_{i+2}}\otimes \cdots\otimes
\ol{a_{r+1}})+f(\ol{a_1}\otimes \cdots \otimes \ol{a_r})\cdot a_{r+1}.
\]
%There is also a normalized version
%$\overline{C}^*(A)=\mathrm{Hom}_{A^e}(\mathbb{B}_*(A), A)\simeq
%\mathrm{Hom}_k(\overline{A}^{\otimes *}, A)$.

For bimodules $M$ and $N$, given $\alpha\in C^p(A, M)$ and
$\beta\in C^q(A, N)$,  the cup product $$\alpha\cup\beta \in
C^{p+q}(A, M\ot_A N)=\mathrm{Hom}_k(\ol{A}^{\otimes (p+q)}, M\ot_A N)$$
  is given by
\[
(\alpha\cup\beta)(\ol{a_1}\otimes \cdots\otimes
\ol{a_{p+q}}):=(-1)^{pq}\alpha(\ol{a_1}\otimes \cdots\otimes
\ol{a_{p}})\ot_A
\beta(\ol{a_{p+1}}\otimes \cdots\otimes
\ol{a_{p+q}}).
\]
This cup product induces a well-defined product in Hochschild
cohomology
\[
\cup \colon HH^p(A, M) \times HH^q(A, N) \longrightarrow
HH^{p+q}(A, M\ot_A N)
\]
which turns the graded $k$-vector space $HH^*(A)=\bigoplus_{n\geq
0}HH^n(A)$ into a graded commutative algebra (\cite[Corollary
1]{Gerstenhaber}).

The Lie bracket is defined as follows. Let $\alpha \in C^n(A, A)$
and $\beta \in C^m(A, A)$. If $n, m\geq 1$, then for $1\leq i\leq
n$, set $\alpha\overline{\circ}_i \beta \in C^{n+m-1}(A, A)$ by
\[
(\alpha\overline{\circ}_i \beta)(a_1\otimes \cdots\otimes
a_{n+m-1}):=\alpha(a_1\otimes \cdots \otimes a_{i-1}\otimes
\beta(a_i\otimes \cdots \otimes a_{i+m-1})\otimes a_{i+m}\otimes
\cdots \otimes a_{n+m-1});
\]
if $ n\geq 1$ and $m=0$, then $\beta\in A$ and for $1\leq i\leq
n$, set
\[
(\alpha\overline{\circ}_i\beta)(\ol{a_1}\otimes \cdots\otimes
\ol{a_{n-1}}):=\alpha(\ol{a_1}\otimes \cdots \otimes \ol{a_{i-1}}\otimes \ol{\beta}
\otimes \ol{a_{i}}\otimes \cdots \otimes \ol{a_{n-1}});
\]
for any  other case, set $\alpha\overline{\circ}_i \beta$ to be
zero. Now define
\[
\alpha\overline{\circ} \beta
:=\sum_{i=1}^n(-1)^{(m-1)(i-1)}\alpha\overline{\circ}_i\beta
\]
and
\[
[\alpha,\, \beta] :=\alpha\overline{\circ}
\beta-(-1)^{(n-1)(m-1)}\beta\overline{\circ} \alpha.
\]
Note that $[\alpha,\, \beta]\in C^{n+m-1}(A, A)$. The above $[\
\,,\,\ ]$ induces a well-defined Lie bracket in Hochschild
cohomology
\[
[\ \,,\,\ ] \colon HH^n(A) \times HH^m(A) \longrightarrow
HH^{n+m-1}(A)
\]
such that $(HH^*(A),\, \cup,\, [\ \,,\,\ ])$ is a Gerstenhaber
algebra (\cite[Page 267, Theorem]{Gerstenhaber}).

The complex used to compute the Hochschild homology $HH_*(A, M)$
is $C_*(A, M)=M\otimes_{A^e} \mathrm{Bar}_*(A)$. Notice that for $r\geq 0$,
$C_r(A, M)=M\otimes_{A^e}(A\ot \ol{A}^{\otimes r}\ot A)\simeq  M\ot \ol{A}^{\otimes r}$
and the differential $$b_r: C_r(A, M)= M\ot \ol{A}^{\otimes
r}\to C_{r-1}(A, M)=M\ot \ol{A}^{\otimes (r-1)}$$  sends $x\otimes
\ol{a_1}\ot  \cdots \otimes \ol{a_{r}}$ to $$\begin{array}{l}  xa_0\ot \ol{a_1}\ot \cdots \ot
\ol{a_r}+\sum_{i=1}^{r-1} (-1)^{i} x\ot \ol{a_1}\otimes \cdots\otimes
\ol{a_{i-1}}\otimes \ol{a_ia_{i+1}}\otimes \ol{a_{i+2}} \otimes \cdots\otimes
\ol{a_r}\\ \ \ \ \  + (-1)^r a_r x\otimes \ol{a_1}\otimes \cdots\otimes \ol{a_{r-1}}.\end{array}$$

There is a  A. Connes' $\mathrm{B}$-operator in the Hochschild
homology theory which is defined as follows. For $a_0\otimes \ol{a_1}\ot
\cdots \otimes \ol{a_r}\in C_r(A, A)$, let $\mathrm{B}(a_0\otimes \ol{a_1}\ot
\cdots \otimes \ol{a_r}) \in C_{r+1}(A, A)$ be
\[
\sum_{i=0}^r (-1)^{ir} 1\otimes \ol{a_i}\otimes\cdots \otimes
\ol{a_r}\otimes \ol{a_0}\otimes \cdots \otimes \ol{a_{i-1}}.
\]
It is easy to check that $\mathrm{B}$ is a chain map satisfying
$\mathrm{B} \comp \mathrm{B}=0$, which induces an operator
$\mathrm{B}: HH_r(A)\rightarrow HH_{r+1}(A)$.

There is a pairing between the Hochschild cohomology and
Hochschild homology, which is called the cap product. For
bimodules $M$ and $N$,  there is a bilinear map
$$\cap: C_r(A, N)\otimes C^p(A, M)\to C_{r-p}(A, N\otimes_A M)$$
sending $z\ot \alpha$ to $$z\cap \alpha=(-1)^{rp} (x\ot_A
\alpha(\ol{a_1}\ot\cdots \ot \ol{a_p}))\ot \ol{a_{p+1}}\ot \cdots \ot \ol{a_r}$$ for
$z=x\ot \ol{a_1}\ot \cdots \ot \ol{a_r}\in C_r(A, N)$ and $\alpha\in C^p(A,
M)$. On verifies easily that $\cap$ induces a well-defined map on
the level of homology, still denoted by $\cap$,
$$\cap: HH_r(A, N)\otimes HH^p(A, M)\to HH_{r-p}(A, N\otimes_A M).$$

I.M.~Gelfand, Yu.L.~Daletskii and B.L.~Tsygan  proved the
following result; see also \cite{TamarkinTsygan}.

\begin{Thm}[\cite{GelfandDaleskiiTsygan}] \label{HochschildCalculus}
The data $(HH^*(A, A), \ \cup, \ [\ ,\ ], \ HH_*(A, A), \ \cap,
\ \mathrm{B})$ is a differential calculus.
\end{Thm}

\subsection{Batalin-Vilkovisky algebras}

In the last  decade, a new structure in Hochschild theory has been
observed, this is the so called Batalin-Vilkovisky structure.

\begin{Def}\label{Def-BV-algebra} A {\em Batalin--Vilkovisky  algebra} (BV algebra for short) is a
Gerstenhaber algebra $(\mathcal{H}^*,\, \cup,\, [\ \,,\,\ ])$
together with an operator $\Delta\colon \mathcal{H}^* \rightarrow
\mathcal{H}^{*-1}$ of degree $-1$ such that $\Delta\comp \Delta=0$, $\Delta(1)=0$
and
\[
[\alpha, \beta]=(-1)^{|\alpha|+1}(\Delta(\alpha\cup \beta)-
\Delta(\alpha)\cup  \beta-(-1)^{|\alpha|}\alpha\cup  \Delta
(\beta)),
\]
for   homogeneous  elements $\alpha, \beta\in \mathcal{H}^*$.
The BV-operator $\Delta:{\mathcal H}^*\ra{\mathcal H}^{*-1}$ is called a
generator of the Gerstenhaber
bracket $[\;,\;]$.
\end{Def}

For a Calabi-Yau algebra, V.~Ginzburg showed in \cite{Ginzburg} that the Hochschild
cohomology of a Calabi-Yau algebra has a Batalin-Vilkovisky
structure. More precisely, for  a
Calabi-Yau algebra $A$ of global dimension $d$, there is a duality
$HH^p(A)\simeq HH_{d-p}(A))$ for $p\geq 0$. Via this duality, we
obtain an operator $\Delta\colon HH^{p}(A)\rightarrow HH^{p-1}(A)$
which is the dual of Connes' operator. This is just the operator
$\Delta$  in the Batalin-Vilkovisky structure.

N.~Kowalzig and
U.~Kr\"{a}hmer extended  the result of V.~Ginzburg  to twisted Calabi-Yau
algebras under a certain condition Let $A$ be a twisted Calabi-Yau algebra
with semisimple algebra automorphism $\sigma$. Then
the Hochschild cohomology ring of $A$ is   a Batalin-Vilkovisky algebra; see \cite[Theorem 1.7]{KowalzigKrahmer}.

T.~Tradler showed  that the Hochschild cohomology algebra of a
symmetric algebra is a BV algebra \cite{Tradler}, see
also~\cite{EuSchedler}\cite{Menichi}. For a symmetric  algebra
$A$, he showed that the $\Delta$-operator on the Hochschild
cohomology corresponds to the Connes' $\mathrm{B}$-operator on
the Hochschild homology via the duality between the Hochschild
cohomology and the Hochschild homology.

\subsection{Algebras with duality, Tamarkin-Tsygan calculi and BV-structures}

Generalising \cite{Lambre} we define.

\begin{Def}
An algebra with duality is given by
$({\mathcal H}^*,  \cup, {\mathcal H}_*,  c,\partial)$, where
\begin{itemize}
\item $({\mathcal H}^*,\cup)$ is a graded commutative unitary
algebra with unit $1\in{\mathcal H}^0$,
\item ${\mathcal H}_*$ is a graded vector space and $c$ is an element
of ${\mathcal H}_d$,
\item $\partial $ is an isomorphism of vector spaces $\partial:{\mathcal H}_*\ra{\mathcal H}^{d-*}$
satisfying $\partial(c)=1$.
\end{itemize}
\end{Def}

Inspired by the result of V.~Ginzburg, the first author gave the
following result which shows for an algebra with duality there is an equivalence
between BV-structure and Tamarkin-Tsygan calculus.

\begin{Lem}\label{Ginzburgrelationlemma}
Let $({\mathcal H}^*,\cup,{\mathcal H}_*, c,\partial)$ be an algebra with duality.
\begin{enumerate}
\item We suppose that
\begin{enumerate}
\item $({\mathcal H}^*,\cup,[\;,\;],{\mathcal H}_*, \cap, B)$ is a Tamarkin-Tsygan calculus,
\item the duality $\partial$ is a homomorphism of ${\mathcal H}^*$-right modules,
i.e. we have the relation
$$\partial(z\cap\alpha)=\partial(z)\cup\alpha.$$
\end{enumerate}
Then the Gerstenhaber algebra $({\mathcal H}^*,\cup,[\;,\;])$ is a BV-algebra
with generator $\Delta=\partial\circ B\circ\partial^{-1}$.
\item We suppose that  $({\mathcal H}^*,\cup,[\;,\;],\Delta)$ is a BV-algebra
with generator $\Delta$.
Then posing $B:=\partial^{-1}\circ \Delta\circ \partial$ and
$z\cap\alpha:=\partial^{-1}(\partial(z)\cup\alpha)$, the data
$({\mathcal H}^*,\cup,[\;,\;],{\mathcal H}_*, \cap, B)$ is a Tamarkin-Tsygan calculus.
\end{enumerate}
\end{Lem}

\begin{Proof}
For (1) see \cite[Lemme 1.6]{Lambre} and (2) is an easy verification.
A similar idea also appeared in \cite[Remark 2.3.67]{EuSchedler}.
\end{Proof}

\begin{Rem}\rm
\begin{enumerate}
\item Regardless its simplicity the relation $\partial(z\cap\alpha)=\partial(z)\cup\alpha,$
which was first noted by V. Ginzburg in \cite[Theorem 34.3]{Ginzburg}, is necessary.
For this reason we call it the Ginzburg relation.
\item This lemma allows to establish the results of V. Ginzburg (for Calabi-Yau algebras), and of
Kowalzig and Kr\"ahmer (for twisted Calabi-Yau algebras).  We shall see that it applies also to
the case of Frobenius algebras.
\end{enumerate}
\end{Rem}

\section{Tamarkin-Tsygan calculus associated with an automorphism of an algebra}

\label{Sectiontwo}

Let $A$ be an associative, finite dimensional and unitary $k$-algebra and let ${\mathfrak N}:A\ra A$ be an
automorphism of this algebra. The aim of this paragraph
is to construct a Tamarkin-Tsygan calculus associated to
this automorphism $\mathfrak N$. Denote by $A_{\mathfrak N}$
the $A$-$A$-bimodule which is $A$ as a $k$-vector space,
and on which we define the bimodule action as $a\cdot m\cdot b:=am{\mathfrak N}(b)$.
Kowalzig and Kr\"ahmer define in \cite[2.18, 7.2]{KowalzigKrahmer} a morphism of $k$-vector spaces
$$\beta_{\mathfrak N}:C_p(A,A_{\mathfrak N})\ra C_{p+1}(A,A_{\mathfrak N})$$
by
$$\beta_\mathfrak{N} (a_0\ot a_1\ot \cdots\ot a_p)=
\Sigma_{i=0}^p(-1)^{ip}\
1\otimes a_{i}\otimes\cdots\otimes a_p
\otimes a_0\otimes \mathfrak{N} (a_1)\otimes\cdots \otimes\mathfrak{N} (a_{i-1})$$
Let $T:C_p(A,A_\mathfrak{N} )\to C_p(A,A_\mathfrak{N} )$ be the morphism defined by
$$T(a_0\ot\cdots\ot a_p)=\mathfrak{N} (a_0)\ot \mathfrak{N} (a_1)\ot \cdots \ot \mathfrak{N} (a_p).$$

\begin{Prop} (N. Kowalzig and U. Kr\"ahmer \cite{KowalzigKrahmer}).
On the space  $C_p(A,A_\mathfrak{N} )$, we get the identity
$$b\beta_\mathfrak{N} +\beta_\mathfrak{N}  b=1-T$$
where $b$ is the Hochschild differential.
\end{Prop}

\begin{Proof} See \cite[2.19]{KowalzigKrahmer} in the setup of Hopf
algebroids; for a proof in the setup of  Hochschild cohomology,
see \cite[Section 4]{GoodmanKraehmer}. \end{Proof}

\subsection{Decomposition of the homology associated with the spectrum of an automorphism}

Let $\Lambda$ be the set of eigenvalues of the automorphism $\mathfrak N$.
Suppose that $\Lambda\subset k$.
Fix an eigenvalue $\lambda \in \Lambda$ of $\mathfrak{N} $ and let $A_\lambda$
be the eigenspace associated
with $\lambda$.
It is trivial to see that for $\lambda, \mu\in \Lambda$,  we get $A_\lambda\cdot A_\mu\subseteq A_{\lambda\mu}$.
When $\lambda\mu\not \in \Lambda$, it is understood that $A_{\lambda\mu}=0$.

For $\lambda\in \Lambda$, write $\overline{A}_\lambda=A_\lambda$ for $\lambda\neq 1$
and $\overline{A}_1=A_1/(k\cdot 1_A)$, and
put $$C_p^\lambda (A,A_\mathfrak{N} ):=\oplus_{\lambda_{i}\in\Lambda, \prod {\lambda_i}=\lambda}
\ \ A_{\lambda_ 0}\otimes \overline{A}_{\lambda_ 1}\otimes\cdots\otimes \overline{A}_{\lambda_ p}.$$
The  Hochschild differential $b:C_p(A,A_\mathfrak{N} )\to C_{p-1}(A,A_\mathfrak{N} )$
restricts to this subspace and denote its restriction by $b^{\lambda}$, then
  $(C_*^\lambda(A,A_\mathfrak{N}), b^\lambda)$
 is a sub-complex of $(C_*(A,A), b) $.
Denote $$HH_{p}^{\lambda} (A,A_\mathfrak{N} ):=H_p(C_*^\lambda (A,A_\mathfrak{N} ), b^\lambda).$$
 We hence obtain a vector space homomorphism
$HH_*^\lambda(A, A_\mathfrak{N} )\to HH_*(A,A_\mathfrak{N} ).$

\begin{Prop}
\begin{enumerate}
\item For each eigenvalue $\lambda\not=1$ of the automorphism $\mathfrak{N} $ we get
$$HH_*^\lambda(A, A_\mathfrak{N} )=0.$$
\item The restriction  $\beta_\mathfrak{N} ^1:C_*^1(A,A_\mathfrak{N} )\to C_{*+1}^1(A,A_\mathfrak{N} )$
of the map $\beta_\mathfrak{N} $ to the sub-complex associated to the eigenvalue $1$
induces a Connes operator
$$B_\mathfrak{N} :HH^1_*(A,A_\mathfrak{N} )\to HH^1_{*+1}(A,A_\mathfrak{N} )$$
with coefficients in the twisted bimodule $A_\mathfrak{N} $, and this map
satisfies $B_\mathfrak{N} ^2=0$.
\end{enumerate}
\end{Prop}

\begin{Proof}
For each eigenvalue $\lambda$ of the automorphism $\mathfrak{N} $  of the algebra $A$, on
$C_p^\lambda(A,A_\mathfrak{N} )$, we obtain the identity
$$b^\lambda \beta_\mathfrak{N} +\beta_\mathfrak{N}  b^\lambda=1-T.$$
However, the restriction of $T$ to $C_*^\lambda(A, A_\mathfrak{N})$ is
$\lambda\cdot \text{ id}$, we get
$b^\lambda \beta_\mathfrak{N} +\beta_\mathfrak{N}  b^\lambda=(1-\lambda) \cdot \mathrm{id}.$
Whenever $\lambda\not=1$ the complex $(C_*^\lambda(A,A_\mathfrak{N} ) , b^\lambda)$
is acyclic with  contracting homotopy $\beta_\mathfrak{N}$.
For $\lambda=1$ we get
$b_1 \beta_\mathfrak{N} ^1+\beta_\mathfrak{N} ^1 b_1=0$, which defines $B_\mathfrak{N}$.
The relation $B_\mathfrak{N} ^2=0$ is a consequence of \cite[2.19]{KowalzigKrahmer}.
\end{Proof}

An analogous decomposition exists for cohomology. Let
$C^p_\lambda (A)$ be those Hochschild cochains $\varphi\in C^p(A,A)$ such that
we have $\varphi( \ol{A}_{\mu_1}\ot \cdots \ot \ol{A}_{\mu_p})\subset A_{\lambda\mu_1\cdots \mu_p}$
for all eigenvalues $\mu_1, \cdots, \mu_p$ of $\mathfrak N$.
The restriction $b_\lambda$ of the Hochschild differential $b:C^p(A,A)\to C^{p+1}(A,A)$ to $C^p_\lambda(A)$
has values in $C^{p+1}_\lambda(A)$. Put
$$HH^p_\lambda(A, A):=H^p(C^*_\lambda(A), b_\lambda).$$
The sub-complex $(C^*_\lambda (A,A),b_\lambda)$ of $(C^*(A,A), b)$ defines a morphism of vector spaces
$HH^*_\lambda(A, A)\to HH_*(A, A)$.

If $\lambda$ and $\mu$ are two eigenvectors of $\mathfrak{N} $ we verify that the cup-product
$\cup:HH^p(A,A)\otimes  HH^q(A,A)\to HH^{p+q}(A,A)$ and the Gerstenhaber bracket
$[\ ,\ ]: HH^p(A,A)\otimes  HH^q(A,A)\to HH^{p+q-1}(A,A)$ induce restrictions
$$\cup_{\lambda,\mu}:HH^p_\lambda(A,A)\otimes HH^q_\mu(A,A)\to HH^{p+q}_{\lambda\mu}(A,A)$$ and
$$[\ ,\ ]_{\lambda,\mu}:HH^p_\lambda(A,A)\otimes HH^q_\mu(A,A)\to HH^{p+q-1}_{\lambda\mu}(A,A).$$
Analogously, the {\em cap}-product
$\cap:HH_p(A,A_\mathfrak{N} )\otimes HH^q(A,A)\to HH_{p-q}(A,A_\mathfrak{N} )$
induce restrictions
$$\cap_{\lambda,\mu}:HH_p^\lambda(A, A_\mathfrak{N} )\otimes
HH^q_\mu(A,A)\to HH_{p-q}^{\lambda\mu}(A,A_\mathfrak{N} ).$$

\subsection{The case of eigenvalue $1$}

Apply the results above to the case $\lambda=\mu=1$. We then get

\begin{Thm}\label{Thm:IntermediateCalculus}
Let $\mathfrak{N} $ be an automorphism of the algebra $A$.
Let $\Lambda$ be the set of eigenvalues of the automorphism $\mathfrak N$.
Suppose that $\Lambda\subset k$. Let $\cup_1:=\cup_{1,1}$,
$[\ ,\ ]_1:=[\ ,\ ]_{1,1}$ and $\cap_1:=\cap_{1,1}$ be the restrictions of the
cup-products, Gerstenhaber bracket and  {\em cap}-product to the homology and
cohomology spaces associated with the eigenvalue $1$.
Then Connes' operator $B_\mathfrak{N} $ gives
$$(HH^*_1(A,A), \cup_1, [\ ,\ ]_1, HH^1_*(A,A_\mathfrak{N} ),\cap_1, B_\mathfrak{N} )$$ the structure of a
Tamarkin-Tsygan calculus.
\end{Thm}

\begin{Rem}\rm
This Tamartin-Tsygan calculus applies in diverse types of
algebras for which its Hochschild cohomology/homology are  naturally equipped
with a duality:
\begin{itemize}
\item
 Calabi-Yau algebras for which for which the dualising module $\cal D$ is isomorphic
to the module $A$. In this case the automorphism
$\mathfrak{N} $ is the identity and this is the situation studied by V.~Ginzburg.

\item Twisted    Calabi-Yau algebras for which the dualising module $\cal D$ is isomorphic
to the module $A_\mathfrak{N}$. This is the situation studied by N. Kowalzig and U. Kr\"ahmer.
\item Symmetric algebras for which  the Nakayama automorphism is
$\mathfrak{N}=\text{id}$. This is the situation studied by T.~Tradler.
\item   Frobenius algebras. This is the situation studied in this paper.
\end{itemize}
\end{Rem}

\subsection{The diagonalisable case}

\begin{Prop}\label{hochschldcohomologyfordiagonalisable}
If $\mathfrak{N}$ is diagonalisable, then
$$HH_*^1(A,A_\mathfrak{N})\simeq HH_*(A,A_\mathfrak{N}).$$
\end{Prop}

\begin{Proof} Since $A=\oplus_{\lambda\in\Lambda}A_\lambda$, we get
$$(C_*(A,A_\mathfrak{N}),b)=\oplus_{\lambda\in\Lambda}(C_*(A,A_\mathfrak{N}),b^\lambda)$$
and therefore
$HH_*(A,A_\mathfrak{N} )=\oplus_{\lambda\in\Lambda}HH_*^\lambda(A,A_\mathfrak{N} )$. For $\lambda\not=1$,
we get $HH_*^\lambda(A,A_\mathfrak{N})=0$. This proves
$HH_*(A,A_\mathfrak{N} )=HH_*^1(A,A_\mathfrak{N} ).$
\end{Proof}

\section{The Hochschild cohomology ring of a Frobenius algebra}

\label{Sectionthree}

\subsection{Algebra with duality associated with a Frobenius algebra}
\label{Sectionhreepoint1}

Let $k$ be a field and let $A$ be a finite dimensional $k$-algebra.
Recall (cf e.g. \cite[Section 1.10.1]{reptheobuch} or \cite{YS})
that $A$ is a Frobenius algebra, if there is a
non-degenerate associative bilinear form $\langle-, -\rangle:
A\times  A\to k$. Here the associativity means that $\langle ab,
c\rangle=\langle a, bc\rangle$ for all $a, b$ and $c$ in $A$. Endow
$D(A)=\mathrm{Hom}_k(A, k)$, the $k$-dual of $A$,  with the
canonical bimodule structure
$$(afb)(c)=f(bca),\ \mathrm{for}\  f\in D(A), a, b, c\in A.$$
The property of being Frobenius is   equivalent to saying  that
$D(A)=\mathrm{Hom}_k(A, k)$ is isomorphic to $A$ as left or as right
modules.  It is readily seen that the map  $a\mapsto \langle a,
-\rangle$ for $a\in A$ gives an   isomorphism of right modules
between    $ A$ and $D(A)$, while the map  $a\mapsto \langle -,
a\rangle$ gives the isomorphism of left modules.   For $a\in A$,
  there exists a unique
$\mathfrak{N}(a)\in A$ such that $\langle a,\  -\rangle=\langle
-,\ \mathfrak{N}(a)\rangle\in D(A)$. It is easy to see that $\mathfrak{N}: A\to
A$ is an algebra isomorphism and we call it
  the Nakayama automorphism of $A$ (associated to
the bilinear form $\langle-, -\rangle$). As above we write $A_\mathfrak{N}$ for the
$A$-$A$-bimodule whose underlying space is $A$ and where the left
$A$-module structure is given by left multiplication and the right
$A$-module structure is given by $x\cdot a=x\ \mathfrak{N}(a)$ for $x\in
A_\mathfrak{N}$ and $a\in A$. Then the map $a\mapsto \langle -, a\rangle$
is an isomorphism of bimodules $ A_\mathfrak{N}\simeq D(A)$. In fact for $x\in
A_\mathfrak{N}$ and $a\in A$, via the isomorphism of left
modules $A\simeq D(A)$,  $x\ \mathfrak{N}(a)$ is sent to
$$     \langle-,\  x\ \mathfrak{N}(a)\rangle=\langle
\mathfrak{N}^{-1}(x\mathfrak{N}(a)),\  -\rangle
=  \langle \mathfrak{N}^{-1}(x)\  a,\   -\rangle= \langle \mathfrak{N}^{-1}(x),\  a\ -\rangle
=\langle a\ -,\  x\rangle=\langle -,\  x\rangle\
a.$$

Using the isomorphism of bimodules $D(A) \simeq  A_\mathfrak{N}$, we can
establish a well known duality between Hochschild cohomology
and Hochschild homology groups. In fact there are isomorphisms of complexes: $$\begin{array}{rl}D(C_*(A,A_\mathfrak{N}))&=
\mathrm{Hom}_k(A_\mathfrak{N}\ot_{A^e}\mathrm{Bar}_*(A), k)\simeq \mathrm{Hom}_{A^e}
(\mathrm{Bar}_*(A), D(A_\mathfrak{N}))\\ &\simeq \mathrm{Hom}_{A^e} (\mathrm{Bar}_*(A), A)=C^*(A,
A),\end{array}$$ where the third isomorphism is induced by the isomorphism $ A_\mathfrak{N}\simeq D(A)$.
This induces an isomorphism
$$\partial:D(HH_*(A, A_\mathfrak{N}))\stackrel{\simeq}{\lra} HH^*(A, A).$$
This isomorphism comes from the pairing $$HH_*(A,A_{\mathfrak N})\otimes HH^*(A,A)\ra k.$$
Explicitly for $a_0\ot \ol{a_1}\ot \cdots \ot \ol{a_p}\in C_p(A, A_\mathfrak{N})$ and
$\alpha\in C^p(A, A)$, the pairing is given by
$$\langle  \ a_0\ot \ol{a_1}\ot \cdots \ot \ol{a_p},
\alpha\rangle= (-1)^{p} \langle a_0, \alpha(\ol{a_1}\ot \cdots\ot \ol{a_p})\rangle.$$

\begin{Rem}\label{frobeniusfundamentalclassremark}\rm
The isomorphism $\partial$ is not easy to describe, but its inverse $\partial^{-1}:HH^*(A,A)\stackrel{\simeq}{\lra}D(HH_*(A,A_{\mathfrak N}))$ is given by
$\partial^{-1}(\alpha)=(-1)^{|\alpha|}\langle-,\alpha\rangle$. In particular for $\alpha=1_A\in HH^0(A,A)$
we put $c:=\partial^{-1}(1_A)=\langle-,1_A\rangle$. In other words, the class $c\in D(HH_*(A,A_{\mathfrak N}))$
is chosen such that $\partial(c)=1_A$.
\end{Rem}

\begin{Def}
The element $c\in D(HH_*(A,A_{\mathfrak N}))$ from Remark~\ref{frobeniusfundamentalclassremark}
is called the fundamental class of the Frobenius algebra $A$.
\end{Def}

\begin{Prop}
Let $A$ be a Frobenius algebra with Nakayama automorphism $\mathfrak N$. Put
${\mathcal H}_{-*}:=D(HH_*(A,A_{\mathfrak N}))$, ${\mathcal H}^*=HH^*(A,A)$ and
$c=\langle-,1_A\rangle\in{\mathcal H}_0$.
\begin{enumerate}
\item There is a cap product $\cap:{\mathcal H}_{-p}\otimes{\mathcal H}^q\ra{\mathcal H}_{-(p+q)}$
for which the isomorphism $\partial^{-1}:{\mathcal H}^*\to{\mathcal H}_*$ is the cap product by the fundamental class, i.e. for all $\alpha\in{\mathcal H}^*$ it satisfies the equation $\partial^{-1}(\alpha)=c\cap\alpha$.
\item The inverse isomorphism $\partial:{\mathcal H}_{-*}\ra{\mathcal H}^*$
is a morphism of ${\mathcal H}^*$-modules
i.e. it satisfies the Ginzburg relation $\partial(z\cap \alpha)=\partial(z)\cup\alpha$.
\end{enumerate}
\end{Prop}

\begin{Proof}
(1). For $z\in {\cal H}_{-p}$ and $\alpha\in {\cal H}^q$, define
$z\cap \alpha\in{\cal H}_{-(p+q)}$ as follows. For $t\in HH_{p+q}(A,A_\mathfrak{N})$
we have  $t\cap\alpha\in HH_p(A,A_\mathfrak{N})$. The map
$z\cap\alpha: HH_{p+q}(A,A_\mathfrak{N})\to k$ is defined by
$(z\cap\alpha)(t):=(-1)^{(p+q)q} z(t\cap\alpha)$, that is,
$$z\cap\alpha(-)=(-1)^{(|z|+|\alpha|)\cdot|\alpha|}z(- \cap \alpha).$$
We claim that  $\alpha\in{\cal H}^{p}$, the equality
$\partial^{-1}(\alpha)=c\cap\alpha$ holds. In fact,
we know from the previous remark that  $\partial^{-1}(\alpha)=(-1)^{|\alpha|}\langle\ -\ ,\alpha\rangle$.
Suppose that $\alpha=cl(f)$ is the cohomology class of
$f:\ol A^{\otimes p}\to A$ and  $u=cl(a_0\otimes \ol a_1\otimes\cdots\otimes \ol a_{p})$
the homological class of
$a_0\otimes \ol a_1\otimes\cdots\otimes \ol a_{p}\in A_\mathfrak{N}\otimes \ol A^{\otimes p}$.
Then we get
$u  \cap  \alpha=(-1)^{p} a_0f(\ol a_1\otimes\cdots\otimes \ol a_{p})\in A$ and
$$(c\cap\alpha)(u)=(-1)^{(0+p)p}c(u\cap\alpha) =(-1)^p
\langle u \cap\alpha,1_A\rangle=(-1)^p\langle (-1)^{p} a_0f(\ol a_1\otimes\cdots\otimes \ol a_{p}),1_A\rangle;$$
on the other hand,
$$\partial^{-1}(\alpha)(u)=(-1)^p\langle u,\ \alpha\rangle=
(-1)^{p}\langle a_0\otimes\ol a_1\otimes\cdots\otimes \ol a_{p},f\rangle=(-1)^{p}(-1)^p\langle a_0,f(\ol a_1\otimes\cdots\otimes \ol a_{p})\rangle=(c\cap\alpha)(u).$$
This proves that for all $\alpha\in HH^p(A,A)$
%$\langle -,\alpha\rangle=c(- \cap \alpha)$, that is,
one has the equality $\partial^{-1}(\alpha)=c\cap\alpha$ in $D(HH_p(A,A_{\mathfrak N}))$.

 \medskip

(2). For $z\in{\cal H}_{-p}$,  $\alpha\in{\cal H}^r$ and $\beta\in{\cal H}^q$, we verify that in
${\cal H}_{-p+q+r}$,  the equality
$z\cap(\alpha\cup\beta)=(z\cap\alpha)\cap\beta$ holds .
It follows that  $\partial^{-1}$ (and also   $\partial$)  is an isomorphism of  ${\cal H}^*$-modules.
\end{Proof}

An alternative and very short proof can be given by
{\em defining} $\cap$ by the Ginzburg relation.
Then $\partial(c\cap\alpha)=\partial(c)\cup\alpha=1\cup\alpha=\alpha$ and we get 1).
The above proof however
gives a much more detailed clarification of the structures in the
sense that the cap product claimed in
the statement of the proposition is indeed
the standard cap product of Hochschild cohomology.

As a whole we obtained the following.

\begin{Prop}\label{Prop:IntermediateAlgebraWithDuality}
Let $A$ be a Frobenius algebra with Nakayama automorphism $\mathfrak N$.
Put
${\mathcal H}_{-*}:=D(HH_*(A,A_{\mathfrak N}))$, ${\mathcal H}^*=HH^*(A,A)$,
$\partial:D(HH_*(A,A_{\mathfrak N}))\stackrel{\simeq}{\lra}HH^*(A,A)$ and
$c=\langle-,1_A\rangle\in{\mathcal H}_0$. Then
$({\mathcal H}^*,\cup,{\mathcal H}_{-*},c,\partial)$ is an algebra with duality.
\end{Prop}

\subsection{The spectrum of the Nakayama automorphism of a Frobenius algebra}

Let $A$ be  a Frobenius algebra with Nakayama automorphism $\mathfrak{N}$.
Let $\Lambda$ be the set of eigenvalues of $\mathfrak{N}$ (in the algebraic
closure $\ol{k}$ of $k$) considered as a
linear transformation of the finite dimensional $k$-vector space $A$.
Notice that elements of $\Lambda$ are not necessarily in $k$.

Since $\mathfrak{N}$ is an automorphism, $0\not\in\Lambda$;
since $\mathfrak{N}(1_A)=1_A$ with $1_A$ the unit element of $A$, we have
$1\in \Lambda$;  for some eigenvectors $x, y\in A$ with eigenvalues
$\lambda, \mu\in \Lambda$ respectively, we have
$\mathfrak{N}(xy)=\mathfrak{N}(x)\mathfrak{N}(y)=\lambda\mu xy$ and  therefore, if $xy\neq
0$ then $\lambda\mu\in \Lambda$.

\begin{Lem} \label{Lem:PropertiesOfLambda} Let $A$ be a Frobenius $k$-algebra with diagonalisable
Nakayama automorphism $\mathfrak{N}$.  Let $\Lambda$ be the set of
eigenvalues of $\mathfrak{N}$. For $\lambda\in \Lambda$, denote by
$A_\lambda$ the corresponding  eigenspace.

\begin{itemize}\item[(i)]   For $\lambda\in \Lambda$, we have  $\lambda^{-1}\in \Lambda$.

\item[(ii)] The isomorphism of bimodules  $D(A)\simeq A_\mathfrak{N}$ induces an
isomorphism of vector spaces $D(A_\lambda)\simeq A_{\lambda^{-1}}$, for any $\lambda \in \Lambda$.

\end{itemize}
\end{Lem}

\begin{Proof}

(i). Since for $0\neq x\in A_\lambda$, $\langle x, -\rangle\in D(A)$ is not the zero
linear transformation and $A=\oplus_{\mu\in \Lambda} A_\mu$,
there exist $\mu\in \Lambda$ and $y\in A_\mu$ such that   $\langle x, y\rangle\neq 0$. Now
$$\langle
x, y\rangle=\langle y, \mathfrak{N}(x)\rangle=\lambda\langle y,
x\rangle=\lambda \langle x, \mathfrak{N}(y)\rangle=\lambda \langle x, \mu y\rangle=\lambda\mu \langle x,
y\rangle.$$ We see that $\lambda\mu=1$ and $\mu=\lambda^{-1}$.
This proves (i) in case that $\mathfrak{N}$ is diagonalisable.

(ii). In course of the proof of (i), we     showed   that  for $\lambda\in \Lambda $ and
$0\neq x\in A_\lambda$, $\langle -,  x \rangle$ is zero on $A_\mu$
for $\mu\neq \lambda^{-1}$.  This shows that $D(A_\lambda)\subseteq A_{\lambda^{-1}}$.
By exchanging the role of $\lambda$ and $\lambda^{-1}$, we get  that the isomorphism $D(A)\simeq
A_\mathfrak{N}$ induces an isomorphism  $D(A_\lambda)\simeq
A_{\lambda^{-1}}$.

\end{Proof}

\begin{Rem} \rm\label{Rem:LambdaIsNotAGroup} In the spirit of Lemma~\ref{Lem:PropertiesOfLambda},
one intends to think that $\Lambda$ is a group. However, this is not true.
A counterexample is given by the algebra
$$A(\lambda)=k\langle X, Y\rangle/(X^2, Y^2, XY-\lambda YX)$$
with $\lambda\in k-\{0\}$.  A direct computation shows that
$\Lambda=\{1, \lambda, \lambda^{-1}\}$ which is not a group unless
$\lambda=1$, or $\lambda$
is a square or cubic root of $1$, i.e. $(\lambda+1)(\lambda^3-1)=0$.
\end{Rem}

\subsection{BV-structure for Frobenius algebras}

Let $A$ be a Frobenius algebra with Nakayama automorphism $\mathfrak N$.
Let $\Lambda$ be the set of eigenvalues of the automorphism $\mathfrak N$
and suppose that $\Lambda\subset k$.
In Section~\ref{Sectiontwo} we obtained a
Tamarkin-Tsygan calculus
$$(HH_1^*(A,A),\cup_1,[\;,\;]_1,HH_*^1(A,A_{\mathfrak N}),\cap_1,B_{\mathfrak N})$$
associated to the eigenvalue $1$ of $\mathfrak N$.
We have constructed  in Section~\ref{Sectionhreepoint1}
the algebra with duality $$({\mathcal H}^*,\cup,{\mathcal H}_{-*},c,\partial).$$
These two structures give an algebra with duality and a Tamarkin-Tsygan
calculus satisfying the
Ginzburg relation from Lemma~\ref{Ginzburgrelationlemma}.

Let ${\mathcal H}_1^*:=HH_1^*(A,A)$ and ${\mathcal H}^1_{-*}:=D(HH_*^1(A,A_{\mathfrak N}))$.
The transpose of the cap product
$$\cap_1:HH_p^1(A,A_{\mathfrak N})\otimes HH_1^q(A,A)\ra HH^1_{p-q}(A,A_{\mathfrak N})$$
yields a cap-product, still denoted by $\cap_1$,
$$\cap_1:{\mathcal H}^1_{-p}\otimes{\mathcal H}_1^q\ra{\mathcal H}^1_{-(p+q)}.$$
We have $c=\langle-,1_A\rangle\in{\mathcal H}_{-0}^1$ and the
restriction $c\cap_1-:{\mathcal H}_1^p\ra{\mathcal H}_{-p}^1$ of
$c\cap-$ to ${\mathcal H}_1^p$ is the isomorphism
$D(HH_p(A, A_\mathfrak{N}))\simeq HH^p_1(A,A)$. This shows
$({\cal H}^*_1, \ \cup_1,\  {\cal H}^*_1, \ c, \ \partial_1)$ is an algebra with duality.
The transpose of Connes' operator
$B_{\mathfrak N}:HH^1_p(A,A_\mathfrak{N})\to HH^1_{p+1}(A,A_\mathfrak{N})$ induces a map
$B_1: {\cal H}_{-(*+1)}^1\to {\cal H}_{-*}^1$.

\begin{Thm}
Let $A$ be a Frobenius algebra with Nakayama automorphism
$\mathfrak N$. Let $\Lambda$ be the set of eigenvalues of the
automorphism $\mathfrak N$. Suppose that $\Lambda\subset k$.
Let $HH_1^*(A,A)$ be the Hochschild cohomology space associated to the eigenvalue $1$
of the Nakayama automorphism $\mathfrak N$. Then the Gerstenhaber algebra $HH_1^*(A,A)$
is a BV-algebra.
\end{Thm}

\begin{Proof} This is because
the algebra with duality $({\cal H}^*_1, \cup_1, {\cal H}_*^1, c, \partial_1)$
and the Tamarkin-Tsygan calculus $({\cal H}^*_1, \cup_1, [\ ,\ ]_1, {\cal H}_*^1, B_1)$
satisfy the hypotheses
of Lemma~\ref{Ginzburgrelationlemma}.
\end{Proof}

\begin{Cor}\label{diagonalisablefrobenius}
Let $A$ be a Frobenius algebra with Nakayama automorphism $\mathfrak N$.
If $\mathfrak N$ is diagonalisable then the Hochschild cohomology $HH^*(A,A)$
is a BV algebra.
\end{Cor}

\begin{Proof}
If $\mathfrak N$ is diagonalisable  we have seen in Proposition~\ref{hochschldcohomologyfordiagonalisable}
that $HH_1^*(A,A)=HH^*(A,A)$.
\end{Proof}

\section{Proof of the main result}\label{Section:mainresult}

\label{Sectionfour}

Let us recall the  statement of our main result of this paper.
\begin{Thm} Let $A$ be a Frobenius algebra with semisimple
Nakayama automorphism. Then  the Hochschild cohomology ring $HH^*(A)$ of
$A$ is  a Batalin-Vilkovisky algebra.
\end{Thm}

The proof of this theorem occupies the rest of this section.

If the Nakayama automorphism is diagonalisable, this is a consequence
of Corollary~\ref{diagonalisablefrobenius}.

Now suppose that the Nakayama automorphism of a Frobenius algebra is semisimple,
that is, it is diagonalisable over the algebraic closure $\ol{k}$ of $k$.

Let $C=A\ot\ol{k}$.  As is readily verified, $C$ is still a
Frobenius algebra with respect to the induced bilinear form
$$\langle a\ot \lambda, b\ot \mu\rangle=\lambda\mu \langle a, b\rangle,
\ \ a, b\in A, \lambda,\mu\in \ol{k}.$$
Therefore, the Nakayama automorphism of $C$ is $\mathfrak{N}_C=\mathfrak{N}\ot \text{id}_{\ol{k}}$.
We shall write $D_{\ol{k}}=\mathrm{Hom}_{\ol{k}} (-, \ol{k})$.

Notice that
$$ D_{\ol{k}}(C)=\mathrm{Hom}_{\ol{k}}(A\ot_k \ol{k}, \ol{k})
\simeq \mathrm{Hom}_{k}(A, \ol{k})\simeq \mathrm{Hom}_k(A, k)\ot \ol{k}=D(A)\ot \ol{k},$$
where the inverse of  the isomorphism $\mathrm{Hom}_{k}(A, \ol{k})\simeq \mathrm{Hom}_k(A, k)\ot \ol{k}$
is given by  $f\ot \lambda \mapsto   (x\mapsto f(x)\ot \lambda)$ for $f\in  \mathrm{Hom}_k(A, k)$
and $\lambda\in  \ol{k}$. We also have  an isomorphism of bimodules
$C_{\mathfrak{N}_C} \simeq A_{\mathfrak{N}}\ot \ol{k}$.   For the Frobenius
$\ol{k}$-algebra $C$, the isomorphism of bimodules $ D(C) \simeq  C_{\mathfrak{N}_C}$
fits into a commutative diagram
$$\xymatrix{  D_{\ol{k}}(C) \ar[r]_{\simeq} \ar[d]_{\simeq} &
C_{\mathfrak{N}_C}\ar[d]^{\simeq} \\ D(A)\ot \ol{k} \ar[r]_{\simeq}  & A_\mathfrak{N}\ot \ol{k}}$$
where the vertical isomorphisms are explicitly given above.

The diagonalisable case of Theorem~\ref{Thm:MainResult} applies to $C$ and therefore
$HH^*_{\ol{k}}(C)$ is a BV algebra, where $HH^*_{\ol{k}}(C)$ is the
Hochschild cohomology of $C$ considered as  a $\ol{k}$-algebra. Denote by
$\Delta_C$ the BV-operator over $HH^*_{\ol{k}}(C)$.

Let us explain the idea of the proof.
It is true  that  $HH^*_{\ol{k}}(C)\simeq HH^*(A)\ot \ol{k}$ as Gerstenhaber algebras; see Proposition~\ref{Prop:FieldExtensionGerstenhaberAlgebra} below.
In order to show that $HH^*(A)$ is a BV algebra,
if we could  show that the $\Delta_C$-operator sends
$HH^*(A)\ot 1=HH^*(A)$ into itself, then  denote by $\Delta_A$ the restriction
of $\Delta_C$ to   $HH^*(A)$.  Since $\Delta_C$ is $\ol{k}$-linear,
we have $\Delta_C=\Delta_A\ot \ol{k}$.

\begin{Prop}\label{Prop:FieldExtensionGerstenhaberAlgebra} Let $A$ be an
algebra defined over a field $k$. Denote $C=A\ot \ol{k}$.
Then there is an isomorphism of Gerstenhaber algebras
$$HH_{\ol{k}}^*(C)\simeq HH^*(A)\ot \ol{k},$$
where  $HH_{\ol{k}}^*(C)$ is the Hochschild cohomology of $C$ considered as a
$\ol{k}$-algebra and the Gerstenhaber algebra structure on $HH^*(A)\ot\ol{k}$
is defined in Remark~\ref{Rem:ExtensionDefinitionOfGertstenhaberAlgebras}.

\end{Prop}
\begin{Proof}  In fact for each  $p\geq 0$,
$$C^p(C, C)=\mathrm{Hom}_{\ol{k}}((C/\ol k\cdot 1)^{\ot_{\ol{k}}p},C)
\simeq \mathrm{Hom}_{\ol{k}}((A/k\cdot 1)^{\ot p}\ot \ol{k}, C)\simeq \mathrm{Hom}_k((A/k\cdot 1)^{\ot p},
C)\simeq C^p(A, A)\ot \ol{k}.$$ One see easily that this is an isomorphism of complexes.
This induces an isomorphism of graded vector spaces $HH_{\ol{k}}^*(C)\simeq HH^*(A)\ot \ol{k}.$
Moreover, a careful examination on the definition of cup product and Lie bracket
shows that this is also an isomorphism of Gerstenhaber algebras.

\end{Proof}

  The proof of the main result then deduces from the following result.

\begin{Lem} \begin{itemize}

\item[(i)] There is an isomorphism of complexes $C_*(C, C_{\mathfrak{N}_C})\simeq C_*(A,
A_{\mathfrak{N}})\ot \ol{k}$.

\item[(ii)] There is an isomorphism of complexes $D(C_*(C, C_{\mathfrak{N}_C}))\simeq D(C_*(A,
A_{\mathfrak{N}}))\ot \ol{k}$.

\item[(iii)] There is a commutative diagram of isomorphisms of complexes
$$\xymatrix{ D_{\ol{k}}(C_*(C, C_{\mathfrak{N}_C})) \ar[r]^\simeq \ar[d]_{\simeq} &  D(C_*(A, A_\mathfrak{N})\ot \ol{k}\ar[d]^\simeq\\
C^p(C, C) \ar[r]^\simeq & C^p(A, A)\ot \ol{k},}$$
where the horizontal isomorphisms are introduced in (i) and (ii), and the
vertical isomorphisms are (induced by) duality isomorphisms.

\item[(iv)]  For each $p\geq 0$, there is a commutative diagram involving Connes
operators over $C$ and $A$
$$\xymatrix{C_p(C, C_{\mathfrak{N}_C})) \ar[r]^{ \mathrm{B}_p} \ar[d]_{\simeq} &  C_{p+1}(C, C_{\mathfrak{N}_C}))\ar[d]^\simeq\\
C_p(A, A_\mathfrak{N})\ot \ol{k} \ar[r]^{\mathrm{B}_p\ot \text{id}_{\ol{k}}}  & C_{p+1}(A, A_\mathfrak{N})\ot \ol{k},}$$
where the vertical isomorphisms are introduced in (ii).
\end{itemize}

\end{Lem}

 \begin{Proof} (i). For each $p\geq 0$,
$$C_p(C, C_{\mathfrak{N}_C})=C_{\mathfrak{N}_C}\ot_{\ol{k}}
(C/\ol k\cdot 1)^{\ot_{\ol{k}} p}\simeq (A_{\mathfrak{N}}\ot (A/k\cdot 1)^{\ot p})\ot \ol{k}= C_p(A,
A_{\mathfrak{N}})\ot \ol{k}.$$
One then verifies that these isomorphisms commute with the differential.

(ii).  For each $p\geq 0$,
$$ D_{\ol{k}}(C_p(C, C_{\mathfrak{N}_C}))=\mathrm{Hom}_{\ol{k}}(C_p(A,
A_\mathfrak{N})\ot \ol{k}, \ol{k})\simeq \mathrm{Hom}_k(C_p(A, A_\mathfrak{N}), \ol{k})\simeq
DC_p(A, A_\mathfrak{N})\ot \ol{k},$$
where  $D_{\ol{k}}$ denotes the $\ol{k}$-dual $\mathrm{Hom}_{\ol{k}}(-, \ol{k})$.

(iii)(iv) The proof can be done by chasing the diagrams.

 \end{Proof}

Now the theorem follows from the diagrams in (iii)(iv) of the above lemma,
since $\Delta$-operator  and Connes operator  $\mathrm{B}$ are dual to each other.

\section{Examples}

\label{Sectionfive}

\subsection{Frobenius algebras with semisimple
Nakayama automorphisms in terms of quiver with relations}

Let $A=kQ/I$ be a finite dimensional  algebra given by quiver with relations.
As is well known, we can choose a basis $\mathcal{B}$ of $A$ consisting of paths
which also contains a basis for the socle of each indecomposable projective $A$-module.
Suppose now that $A$ is a Frobenius algebra. Then by \cite[Proposition 2.8]{HolmZimmermann},
there is a natural choice of the defining bilinear form $\langle a, b\rangle=tr(ab)$
for $a, b\in A$  induced by the trace map
$$tr: A\to k,\ \ p\in \mathcal{B}\mapsto \left\{\begin{array}{cl} 1 &
\mathrm{if}\ p\in Soc(A)\cap \mathcal{B}\\ 0 & \mathrm{otherwise}\end{array}\right.$$

Assume that the basis $\mathcal{B}$ satisfies two further conditions:
\begin{itemize}\item[(1)]  for arbitrary two paths $p, q\in \mathcal{B}$, there exist
another path $r\in \mathcal{B}$ and a constant $\lambda\in k$ such that
$p\cdot q=\lambda r\in A$
\item[(2)] for each path $p\in \mathcal{B}$, there exists a unique element $p^*\in \mathcal{B}$
such that $0\neq p\cdot p^*\in Soc(A)$
\end{itemize}
We can prove the following rather useful result.

\begin{Crit}\label{Crit} Within the above setup, suppose that  $k$ is a field
of characteristic $0$ or  of characteristic $p$
with $p$ strictly biggar than the number of arrows of $Q$.  Then
the two conditions (1) and
(2) imply that the Nakayama automorphism  of $A$ is semisimple and the
Hochschild cohomology of $A$ is a BV algebra.
\end{Crit}

\begin{Proof}
For $p\in \mathcal{B}$, by (2), let $p^*$ be the unique path in
$\mathcal{B}$ such that $p\cdot p^*=\lambda(p) r\in Soc(A)$
with $\lambda(p)\in k\setminus\{0\}$ and $r\in Soc(A)\cap \mathcal{B}$.
Then  for $p, q\in \mathcal{B}$ we get,
$$ \langle p, q\rangle=\left\{\begin{array}{cl} \lambda(p) & \mathrm{if}\ q=p^*\\ 0 & \mathrm{otherwise}\end{array}\right.$$
Since $\langle p, q\rangle=\langle q, \mathfrak{N}(p)\rangle$,
the Nakayama automorphism sends $p$ to $\frac{\lambda(p)}{\lambda(p^*)}p^{**}$.
Since $\mathcal{B}$ is finite,  the Nakayama automorphism $\mathfrak{N}$,
restricted to $\mathcal{B}$, is a  permutation of $\mathcal{B}$,
modulo scalars.

%Hence, a sufficiently large power of $\mathfrak{N}$ is diagonalisable. Since the base field is
% of characteristic $0$ or  of characteristic $p$
%with $p$ big enough, $\mathfrak{N}$ itself is diagonalisable.

We will show that the Nakayama automorphism $\mathfrak N$ is diagonalisable if $k$ is an
algebraically closed field of characteristic $0$ or algebraically closed of characteristic $p$
where $p>\dim(rad(A)/rad^2(A))=:d$. Recall that the arrows $Q_1$ of the quiver of $A$
form a $k$-basis of $rad(A)/rad^2(A)$.
Since $\mathfrak N$ is an algebra automorphism, and since
$A$ satisfies the conditions (1) and (2), for each $p\in Q_1$ we get $p^{**}\in Q_1$.
We will show that
the action of $\mathfrak N$ on the $k$-vector space $M=rad(A)/rad^2(A)$
generated by $Q_1$ is diagonalisable.
Let $G$ be the infinite cyclic group, generated by $c$. Then $kG$ acts on $rad(A)/rad^2(A)$
when we define the action of $c$ on $M$ by $\mathfrak N$.

Let $\alpha\in Q_1$. Then there is a $t_\alpha\in\N\setminus\{0\}$ such that
$c^{t_{\alpha}}\cdot \alpha=u_\alpha\cdot \alpha$ for some
$u_\alpha\in k\setminus\{0\}$. Choose $t_\alpha$ minimal possible.
Let $x_i:=c^i\cdot \alpha$ for $i\in\{0,1,\dots,t_\alpha-1\}$. The $k$-vector space $T_\alpha$
generated by $x_0,\cdots,x_{t_\alpha-1}$ is then a $kG$-module
and $c$ acts by the matrix $C_\alpha$, say.
Using the basis $\{x_0,\cdots,x_{t_\alpha-1}\}$ of $T_\alpha$
it is easily seen that the characteristic polynomial of $C_\alpha$ is $X^{t_\alpha}\pm u$ and
this polynomial
has only simple roots in $k$ since the characteristic of $k$ is either
$0$ or bigger than $d$ and $d\geq t_\alpha$.
Now $M=\bigoplus_{\text{some }\alpha}T_\alpha$.
Let $Q_1'$ be the basis of $M$ for which the action of $\mathfrak N$ is given by a diagonal matrix.
This shows that $\mathfrak N$ acts diagonally on all paths formed by the elements in $Q_1'$.
We may suppose that $A$ is indecomposable as algebra (i.e. $Q$ is connected)
since the Nakayama automorphism acts on each
indecomposable factor. Let $Q_0$ be the set of vertices in the quiver.
If $A$ is indecomposable, then $|Q_1|\geq |Q_0|-1$ and equality holds if and only if $Q$ is a tree.
The quiver of a selfinjective algebra is not a tree, and hence $|Q_1|\geq |Q_0|$.
Since $\mathfrak N$ permutes $Q_0$, the action of $\mathfrak N$ on $kQ_0$ is diagonalisable, using that
the characteristic of the field is $0$ or bigger than $|Q_1|$.
A basis of $A$ is given by $Q_0$ and paths
of elements of $Q_1$. Let $Q_0'$ be a basis of $kQ_0$ and let $Q_1'$ be a basis of $M$
with diagonal action of $\mathfrak N$. Then $\mathfrak N$ acts diagonally on paths produced by elements
of $Q_1'$ and the set of paths of elements of $Q_1'$ forms a generating set of $rad(A)$. Eliminating
superfluous elements we produce this way a basis ${\mathcal B}_r$
of $rad(A)$ on which $\mathfrak N$ acts diagonally.
Hence ${\mathcal B}_r\cup Q_0'$ is a basis of $A$ on which $\mathfrak N$ acts diagonally.
By our main result
Theorem~\ref{Thm:MainResult} the Hochschild cohomology of $A$ is a BV algebra.
\end{Proof}

These seemly rather strong conditions (1) and (2)  are  in fact satisfied by
many interesting classes of   algebras.

\subsection{Tame Frobenius algebras}

In this subsection $k$ denotes an algebraically closed field.

\begin{Lem} Each  self-injective algebra of
finite representation type is Morita equivalent to an algebra $kQ/I$
given by a quiver $Q$ modulo
admissible relations $I$ verifying  the conditions (1) and (2).
\end{Lem}

\begin{Proof}  Each representation-finite algebra has a
multiplicative basis (cf. \cite{BautistaGabrielRoiterSalmeron}), thus the first condition holds.
For the second condition, suppose that for a path $p\in \mathcal{B}$,
there exist two paths $q_1, q_2\in \mathcal{B}$ such that
$0\neq pq_1=\lambda pq_2\in Soc(A)$ with $\lambda\in k$. We can assume that $p$
has positive length, otherwise $q_1$ and $q_2$ would not be linearly independent in $A$,
using that the socle of each  indecomposable  projective module is one-dimensional.
Now $q_1$ and $q_2$ are parallel paths, by reducing suitably their lengths and
enlarging $p$ if necessary, one
can assume that they have no common arrows.  However, this shows that $A$ is of
infinite representation type, as there are infinitely many string modules of
the form $M((q_1q_2^{-1})^n), n\in \mathbb{N}$, which  is a contradiction.

One can also prove this result using a case-by-case analysis based on the list
given  in terms of quiver with relations in \cite{Asashiba}.
\end{Proof}

However, the Nakayama automorphism of a self-injective algebra of
finite representation type is not necessarily semisimple.

\begin{Ex}\rm  \label{Counterexample} Let $k$ be a field of characteristic two.
Consider the algebra defined by the quiver with relations
$$\xymatrix{1\ar@<1ex>[r]^\alpha & 2\ar@<1ex>[l]^{\beta}  & \alpha\beta\alpha\beta=0=\beta\alpha\beta\alpha}$$
Thus $A$ is a self-injective Nakayama algebra.
Then the indecomposable projective $A$-modules are uniserial and has the following form
$$
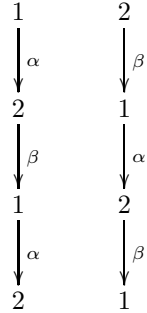
\xymatrix{ 1\ar[d]^\alpha \\ 2\ar[d]^\beta  \\ 1\ar[d]^\alpha  \\  2  }\ \ \ \ \ \ \ \xymatrix{  2\ar[d]^\beta  \\ 1\ar[d]^\alpha  \\  2  \ar[d]^\beta \\  1 }$$
Under the basis $\{e_1, e_2, \alpha, \beta, \alpha\beta,  \beta\alpha\}$,
the matrix of the Nakayama automorphism is
$$\left(\begin{array}{cccccc} 0 & 1 & 0 & 0& 0& 0\\  1 & 0 & 0 & 0& 0& 0\\ 0 & 0 & 0 & 1& 0& 0\\ 0 & 0 & 1 & 0& 0& 0\\ 0 & 0 & 0 & 0& 0& 1\\ 0 & 0 & 0 & 0& 1& 0\end{array}\right).$$
Therefore, the Nakayama automorphism of $A$ is not semisimple.
It would be interesting to see whether $HH^*(A)$ is a BV algebra.

For each prime number $p$, one can construct such a selfinjective Nakayama
algebra  over  field of  characteristic $p$.
\end{Ex}

Another class of algebras is the so-called self-injective special biserial algebras.
A pair $(Q,I)$ of a quiver $Q$ and admissible relations $I$
is called special biserial, if the following conditions hold:
\begin{itemize}
\item[(a)] Each vertex has at most two leaving arrows and at most
two entering arrows.

\item[(b)] Given an arbitrary arrow  $\alpha$, there exists at most one arrow  $\beta$
such that  $t(\alpha)=s(\beta)$ and
$\alpha \beta \not \in I$
and at most one arrow   $\gamma$ such that  $t(\gamma)=s(\alpha)$
and $\gamma \alpha\not \in I$.
\end{itemize}
An algebra is called a special biserial algebra if it is Morita equivalent to
$kQ/I$ for a special biserial pair $(Q,I)$.

\begin{Lem} \label{specialbiserial}
For a special biserial pair $(Q,I)$
the algebra $kQ/I$ satisfies the two conditions (1) and (2).
\end{Lem}

\begin{Proof} It is not difficult to see, and actually well-known
(cf e.g.  \cite{Erdmann}),
that an indecomposable projective module over a self-injective special biserial
algebra is either a uniserial module or a module for which the quotient of the
radical by its socle is the direct sum of two uniserial modules. The first
case is induced by a monomial relation and the second by a commutation  relation.
For the choice of the basis $\mathcal{B}$, one simply  takes representatives of elements in
$kQ/I$ given by paths except that for each indecomposable projective non uniserial
module, where we choose one of  the two paths from its top to its socle. Now the
two conditions hold trivially.
\end{Proof}

Now let us look at weakly symmetric algebras of domestic representation type.
R.~Bocian, T.~Holm and A.~Skowro\'nski
\cite{BocianHolmSkowronski2004,BocianHolmSkowronski2007,HolmSkowronski2006}
classified all weakly symmetric algebras of domestic type over $k$ up to
derived equivalence  and the last two authors of the present paper
gave a classification up to stable equivalences(\cite{ZhouZimmermann}).
In Bocian-Holm-Skowro\'nski classification, a domestic weakly
symmetric standard algebra with singular  Cartan matrix is derived equivalent to
the trivial extension $T(C)$ of a canonical algebra $C$ of Euclidean type
and is thus symmetric;  see \cite[Theorem 1]{BocianHolmSkowronski2004}.
By \cite[Theorem 2]{BocianHolmSkowronski2004} a domestic weakly
symmetric standard algebra with nonsingular  Cartan matrix is derived
equivalent to some  algebras explicitly given in terms of quiver with relations,
denoted by $A(\lambda), A(p, q), \Lambda(n)$ and $\Gamma(n)$.
Note that these algebras are symmetric except $A(\lambda)$ with $\lambda\in k\setminus\{ 0, 1\}$.
However,
$$A(\lambda)=k\langle X, Y\rangle/(X^2, Y^2, XY-\lambda YX)$$
for $\lambda\not\in\{ 0, 1\}$ has a semisimple Nakayama automorphism,
given by a diagonal matrix with coefficients $(1,\lambda^{-1},\lambda,1)$
with respect to the basis $\{1,X,Y,XY\}$
as is easily verified. One may
use the result of the next subsection, as $A(\lambda)$ is a quantum complete
intersection.

By \cite[Theorem 1]{BocianHolmSkowronski2007}\label{DomesticNonStandard}
any nonstandard
self-injective algebra of domestic type  is derived equivalent
(and also  stably equivalent) to an algebra $\Omega(n)$ with
$n\geq 1$.  Let us recall the quiver with relations of $\Omega(n)$.
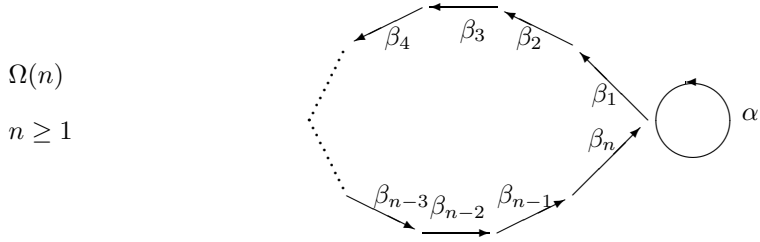
$$\unitlength0.5cm
\begin{picture}(30,11)

\multiput(9, 7)(0.1,0.2){10}{\circle*{0.01}}

\multiput(9, 7)(0.1,-0.2){10}{\circle*{0.01}}

\put(16.5, 7.5){$\beta_1$}

 \put(18, 7){\vector(-1,1){1.8}}

 \put(14.5, 9){$\beta_2$}

\put(16, 9){\vector(-2,1){1.8}}

\put(13, 9.2){$\beta_3$}

\put(14, 10){\vector(-1,0){1.8}}

\put(11, 9.0){$\beta_4$}

\put(12, 10){\vector(-2,-1){1.8}}

\put(10.7, 4.8){$\beta_{n-3}$}

\put(10, 5){\vector(2,-1){1.8}}

\put(12.2, 4.5){$\beta_{n-2}$}

\put(12, 4){\vector(1,0){1.8}}

\put(14, 4.8){$\beta_{n-1}$}

\put(14, 4){\vector(2,1){1.8}}

\put(16.4, 6.3){$\beta_{n}$}

\put(16, 5){\vector(1,1){1.8}}

\put(19.2, 7){\circle{2}}

\put(19, 8){\vector(-1,0){0.01}}

\put(20.5, 7){$\alpha$}

 \put(1, 8){$\Omega(n)$}

 \put(1, 6.5){$n\geq 1$}

 \put(3, 2){$\alpha^2=\alpha\beta_1\beta_2\cdots\beta_n=
 -\beta_1\beta_2\cdots\beta_n\alpha, $}

\put(3, 0.5){$\beta_n\beta_{1}=0, \beta_j\beta_{j+1}\cdots
\beta_n\beta_1\cdots \beta_n \alpha\beta_1 \cdots \beta_{j-1}\beta_j=0,
2\leq j\leq n$}
\end{picture}$$

Notice that we cannot use the Criterion~\ref{Crit} for the algebra
$\Omega(n)$. However, we can still prove the semisimplicity of its Nakayama automorphism.

\begin{Lem} The Nakayama automorphism of $\Omega(n)$ is  diagonalisable.

\end{Lem}

\begin{Proof}  The indecomposable projective modules of $\Omega(n)$ are of the following shape:
$$
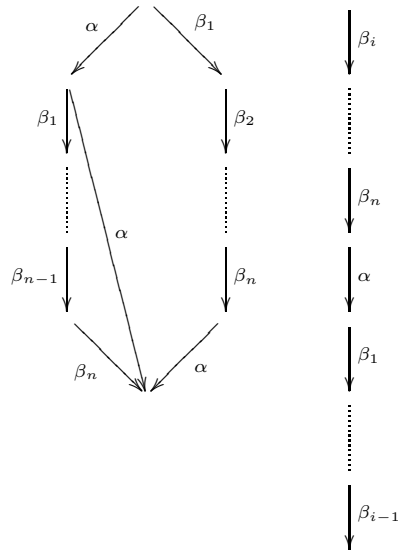
\xymatrix{  & \ar[dl]_\alpha\ar[dr]^{\beta_1} & \\  \ar[d]_{\beta_1}\ar[ddddr]^\alpha & & \ar[d]^{\beta_2} \\ \ar@{..}[d] & & \ar@{..}[d]\\ \ar[d]_{\beta_{n-1}}&& \ar[d]^{\beta_n}\\ \ar[dr]_{\beta_n} && \ar[dl]^\alpha \\ &&  }\ \  \ \ \ \ \ \
\xymatrix{\ar[d]^{\beta_i}\\ \ar@{..}[d] \\   \ar[d]^{\beta_n}   \\  \ar[d]^\alpha \\ \ar[d]^{\beta_1}\\  \ar@{..}[d]\\ \ar[d]^{\beta_{i-1}}  \\ \  }$$
with $2\leq i\leq n$.
This algebra does not satisfy the two conditions of Criterion~\ref{Crit},
but we can compute explicitly its
Nakayama automorphism. For  $\mathcal{B}$, one can take the obvious basis containing
$\alpha, \beta_1\cdots \beta_n, \alpha\beta_1\cdots \beta_n=\alpha^2=-\beta_1\cdots\beta_n\alpha$
etc.  However, the dual basis $\mathcal{B}^*$ does not consist of paths. In fact, one obtains
$\alpha^*=\beta_1\cdots \beta_n$ and $(\beta_1\cdots\beta_n)^*=-\alpha+\beta_1\cdots \beta_n$
etc.  From this, the Nakayama automorphism is given by
$\mathfrak{N}(\alpha)=-\alpha+2\beta_1\cdots \beta_n$ and for any other path
$p\in \mathcal{B}$, we get
$\mathfrak{N}(p)=p$.  Hence, in characteristic two, the Nakayam
automorphism is the identity map (in fact $\Omega(n)$ is symmetric),
and in odd characteristic  it is diagonalisable.  Therefore,
the Nakayama automorphism of $\Omega(n)$ is diagonalisable.
\end{Proof}

We have shown that each weakly symmetric  algebra of domestic type is derived
equivalent to a weakly symmetric  algebra of domestic type whose  Nakayama automorphism
is semisimple.

\bigskip

Now we consider self-injective algebras of polynomial growth which are not of domestic type.
 The derived equivalence classification of the standard (resp.
non-standard) non-domestic weakly symmetric  (resp. self-injective)
algebras of polynomial growth over $k$
is achieved in \cite[Page 653
Theorem]{BialkowskiHolmSkowronski2003a} (resp.\cite[Theorem
3.1]{BialkowskiHolmSkowronski2003b}).

By \cite[Page 653 Theorem]{BialkowskiHolmSkowronski2003a}£¬  an
indecomposable standard non-domestic weakly symmetric algebra of
polynomial growth is always  derived equivalent  to a symmetric algebra
except that it may be derived equivalent to $\Lambda_9'$ in characteristic
not two. For the quiver with relations of the algebra $\Lambda_9'$, we refer
to  \cite{BialkowskiHolmSkowronski2003a}. From this description, we know that
$\Lambda_9'$ is the  preprojective algebra of type $D_4$ and that its Nakayama
automorphism is diagonalisable (and is of order two) by \cite[Section 5.2.1]{Eu};
see also Example~\ref{Ex:Preprojective}.

By  \cite[Theorem 3.1]{BialkowskiHolmSkowronski2003b},
an indecomposable non-standard non-domestic  self-injective
algebra of polynomial growth is always derived equivalent to a
symmetric algebra except the possibility of $\Lambda_{10}$ in characteristic two.
Let us recall its quiver with relation $\Lambda_{10}=kQ/I$:

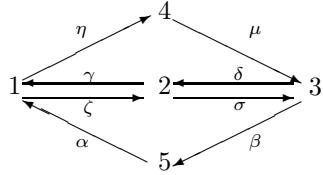
\unitlength1cm

\begin{picture}(10,3)
\put(3,1.5){$1$}
\put(5,1.5){$2$}
\put(5,2.5){$4$}
\put(5,0.5){$5$}
\put(7,1.5){$3$}
\put(3.2,1.7){\vector(2,1){1.7}}
\put(5.2,2.5){\vector(2,-1){1.7}}
\put(6.8,1.65){\vector(-1,0){1.6}}
\put(4.8,1.65){\vector(-1,0){1.6}}
\put(3.2,1.45){\vector(1,0){1.6}}
\put(5.2,1.45){\vector(1,0){1.6}}
\put(6.9,1.4){\vector(-2,-1){1.7}}
\put(4.85,0.6){\vector(-2,1){1.65}}
\put(3.9,2.3){\scriptsize $\eta$}
\put(6.2,2.3){\scriptsize $\mu$}
\put(6.2,0.8){\scriptsize $\beta$}
\put(3.9,0.8){\scriptsize $\alpha$}
\put(4,1.7){\scriptsize $\gamma$}
\put(6,1.7){\scriptsize $\delta$}
\put(6,1.3){\scriptsize $\sigma$}
\put(4,1.25){\scriptsize $\zeta$}

\end{picture}

and
\begin{eqnarray*}
 \Lambda_{10}&=&KQ/(\beta\alpha-\delta\gamma,\zeta\sigma-\eta\mu,\alpha\eta,\mu\beta, \sigma\delta-\gamma\zeta,\delta\sigma\delta\sigma)
\end{eqnarray*}

 \begin{Lem} The Nakayama automorphism of $\Lambda_{10}$ in characteristic two is  not semisimple.

\end{Lem}

\begin{Proof} The indecomposable projective modules of $\Lambda_{10}$ are of the following shape:
$$
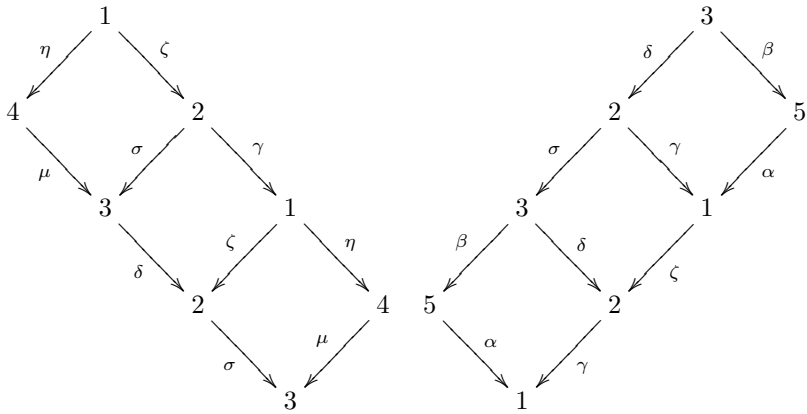
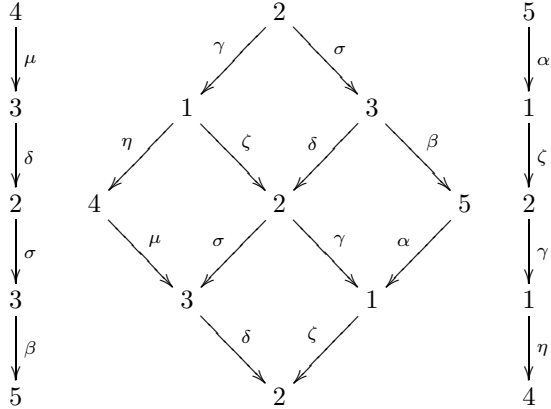
\xymatrix{ & 1\ar[dr]^{\zeta}\ar[dl]_\eta & & & \\ 4 \ar[dr]_\mu & & 2\ar[dr]^\gamma\ar[dl]_\sigma & & \\ & 3\ar[dr]_\delta & & 1\ar[dr]^\eta\ar[dl]_\zeta& \\ & & 2\ar[dr]_\sigma& & 4\ar[dl]_\mu \\ & & & 3 & } \ \
\xymatrix{ & & &  3\ar[dl]_{\delta}\ar[dr]^\beta &  \\ &&  2 \ar[dr]^\gamma\ar[dl]_\sigma & & 5\ar[dl]^\alpha\\
& 3\ar[dr]^\delta \ar[dl]_\beta& & 1\ar[dl]^\zeta&  \\ 5 \ar[dr]^\alpha & & 2\ar[dl]^\gamma & & &  \\ & 1 & & &    } $$ $$
\xymatrix{ 4\ar[d]^\mu \\ 3\ar[d]^\delta  \\ 2\ar[d]^\sigma  \\  3\ar[d]^\beta  \\  5}\ \ \ \
\xymatrix{& & 2\ar[dr]^\sigma\ar[dl]_\gamma & & \\ & 1\ar[dr]^\zeta\ar[dl]_\eta & & 3\ar[dr]^\beta\ar[dl]_\delta & \\ 4 \ar[dr]^\mu& & 2\ar[dr]^\gamma\ar[dl]_\sigma& & 5\ar[dl]_\alpha\\ & 3\ar[dr]^\delta && 1\ar[dl]_\zeta& \\ && 2&&}\ \ \ \
\xymatrix{5\ar[d]^\alpha\\ 1\ar[d]^\zeta \\  2\ar[d]^\gamma  \\  1\ar[d]^\eta  \\ 4}$$

Notice that in the above diagrams, each square is commutative. From this,
we observe that  $\mathfrak N$ is of order $2$.
However, one sees that the Nakayama automorphism permutes the vertices
$1$ and $3$, hence its matrix under a suitable basis has a block
$\left(\begin{array}{cc} 0 & 1\\ 1 & 0\end{array}\right)$ and this matrix
is not diagonalisable in characteristic two.

\end{Proof}

Since derived equivalent algebras have isomorphic Hochschild
cohomology rings (\cite{Ri,Keller}),
we have proved  in this subsection the following

\begin{Prop}
Let $A$ be an algebra falling into  one of the following
classes of algebras
\begin{itemize}
\item  representation-finite   self-injective  algebras in characteristic zero;
\item self-injective  special biserial algebras in characteristic zero;
\item   standard  weakly symmetric algebras  of
domestic type which are not  representation-finite;
\item  nonstandard   self-injective  algebras of
domestic type which are not  representation-finite;
\item   standard non-domestic weakly symmetric algebras  of
polynomial growth;
\item   nonstandard non-domestic self-injective  algebras of
polynomial growth over fields of characteristic different from
$2$, or over fields of characteristic $2$
as long as they are not derived equivalent to $\Lambda_{10}$;
\end{itemize}
Then $A$ is derived equivalent to a Frobenius algebra whose Nakayama
automorphism is semisimple. Therefore,  the Hochschild cohomology ring
of $A$ is a BV algebra.
\end{Prop}

We actually proved a slightly more precise statement concerning the characteristic of $k$.

We do not know whether  BV structures exists or not on the Hochschild
cohomology ring of Example~\ref{Counterexample} or that of $\Lambda_{10}$ in characteristic $2$.

\subsection{Quantum complete intersections}

In \cite{Happel}, D.~Happel asked whether an algebra has finite global dimension
whenever its Hochschild cohomology is finite dimensional. Although Happel's conjecture was
verified for many classes of algebras, it is wrong in general. A counter-example
was exhibited in \cite{BGMS}. This example is in fact our algebra $A(\lambda)$
from Remark~\ref{Rem:LambdaIsNotAGroup}.

This example  has been generalized the so-called \textit{quantum complete intersections},
which are extensively studied by P.A.~Bergh, K.~Erdmann, S.~Oppermann etc. Let $N\geq 2$
and $\textbf{a}=(a_1, \cdots, a_N)$ with $a_j\geq 1$. Let
$\textbf{q}=(q_{ij}, 1\leq i, j\leq N)$ be a family of nonzero constants in $k$ such that
$q_{ii}=1$ and $q_{ij}q_{ji}=1$. Now define
$$A(\textbf{q}, \textbf{a})=\frac{k\langle X_1, \cdots, X_N\rangle}{(X_i^{a_i+1},
X_iX_j-q_{ij} X_jX_i, 1\leq i, j\leq N)}.$$
Obviously this algebra is a local weakly symmetric algebra, and is
thus a Frobenius algebra.  A direct computation
shows that for each $1\leq i\leq N$,
$\mathfrak{N}(X_i)=\big(\prod_{j=1}^N q_{ij}^{a_j}\big) X_i$ and so
it is diagonalisable.

\begin{Cor}
The Hochschild cohomology ring of a quantum  complete intersection
$A(\textbf{q}, \textbf{a})$ is a BV algebra.
\end{Cor}

\subsection{Finite dimensional Hopf algebras}

Let $k$ be an algebraically closed field of characteristic zero. Let $H$ be a
finite dimensional Hopf algebra over $k$.
By \cite{LarsonSweedler} we get that $H$ is Frobenius.
Indeed, given a
right integral $\varphi\in H^*$,  a Frobenius bilinear
form is given by
$\langle a, b\rangle=\varphi(ab)$. Since the antipode $S$ of $H$ has finite order by \cite{Radford},
its Nakayama automorphism  also has finite order.

\begin{Cor}
The Hochschild cohomology ring of a finite dimensional Hopf algebra
defined over an algebraically
closed field of characteristic zero is a Batalin-Vilkovisky algebra.
\end{Cor}

It would  be an interesting question to know  when the usual cohomology
groups of $H$ is a BV subalgebra of $HH^*(H)$;
a sufficient condition was provided by L.~Menichi in \cite[Theorem 50]{Menichi}.

\subsection{Other examples}

There are many other examples of Frobenius algebras related to Calabi-Yau algebras
and Artin-Schelter regular algebras.

\begin{Ex}  \rm In the classical paper \cite{ArtinSchelter}, M. Artin and W.F.~Schelter
classified three dimensional  Artin-Schelter regular algebras. These algebras are
twisted Calabi-Yau algebras, which implies that there is an algebra automorphism
$\sigma$ of $A$ such that $HH^{d-*}(A)\simeq HH_*(A, A_{\sigma})$. In the classification,
they use a generic condition which implies   the semisimplicity of the algebra
automorphisms $\sigma$ of  these algebras. When these algebras are Koszul,
their Koszul duals are Frobenius by
\cite[Corollary D]{LuPalmeriWuZhang} and
the Nakayama automorphism of $A^!$
and the algebra automorphism $\sigma$  of $A$ are related by
\cite[Theorem 9.2]{vandenBergh}. Therefore, whenever $\sigma$ is semisimple,
the Koszul duals are Frobenius algebras with semisimple Nakayama automorphisms.
 The Hochschild cohomology ring of the Koszul dual of a
 three dimensional  Artin-Schelter regular algebra  is  BV algebra.
 We do not know the explicit BV structure over the Hochschild cohomology rings of these algebras.
\end{Ex}

\begin{Ex}\label{Ex:Preprojective} \rm The preprojective algebras of Dynkin
quivers $ADE$ are Frobenius algebras whose Nakayama automorphism are
has finite order; for details see \cite{Eu}.  Except the cases that
$\mathrm{char}\  k=2$, and the type $D_n, n$ odd or $E_6$, the Nakayama
automorphism is diagonalisable.
Therefore, except these cases  their Hochschild cohomology rings are
BV algebras. This is a well known fact (at least over a field of
characteristic zero) and our main result gives a structural
explanation of the existence of BV structure.
This BV structure (over a field of characteristic zero) has
been computed by C.-H.~Eu in \cite{Eu2}.
\end{Ex}

\begin{Ex} \rm
Another class of Frobenius algebras, called almost Calabi-Yau
algebras, was extensively studied by  D.E.~Evans and M.~Pugh
(cf.  \cite{EvansPugh1, EvansPugh2}). These algebras are related to $SU(3)$
modular invariants and MacKay correspondence.  Their Nakayama automorphisms
 have also finite order and is thus semisimple over a field of
 characteristic zero; the authors in fact works over $\mathbb{C}$.
 Therefore,  the Hochschild cohomology ring of an almost Calabi-Yau algebra
 defined over a field of characteristic zero is a BV algebra.
 It would  be interesting to
compute the BV structure over the Hochschild cohomology rings of these algebras.
\end{Ex}

\end{document}